\input amstex\documentstyle {amsppt}  
\pagewidth{12.5 cm}\pageheight{19 cm}\magnification\magstep1
\topmatter
\title Character sheaves and generalizations\endtitle
\author G. Lusztig\endauthor
\address Department of Mathematics, M.I.T., Cambridge, MA 02139\endaddress
\dedicatory{Dedicated to I. M. Gelfand on the occasion of his 90th birthday}
\enddedicatory
\thanks Supported in part by the National Science Foundation\endthanks
\endtopmatter   
\document

\define\Lie{\text{\rm Lie }}

\define\frl{\forall}

\define\si{\sim}

\define\qua{\quad}

\define\tcl{\ti\cl}

\define\op{\oplus}

\define\m{\mapsto}
\define\do{\dots}

\define\sm{\smallmatrix}
\define\esm{\endsmallmatrix}
\define\sub{\subset}
\define\bxt{\boxtimes}
\define\T{\times}
\define\ti{\tilde}
\define\nl{\newline}
\redefine\i{^{-1}}

\define\ot{\otimes}
\define\bbq{\bar{\QQ}_l}

\define\Hom{\text{\rm Hom}}

\define\tr{\text{\rm tr}}

\define\supp{\text{\rm supp}}

\define\a{\alpha}

\redefine\c{\chi}

\define\e{\epsilon}

\define\p{\pi}
\define\ph{\phi}
\define\ps{\psi}
\define\r{\rho}

\define\x{\xi}

\define\kk{\bold k}

\define\FF{\bold F}

\define\QQ{\bold Q}

\define\ZZ{\bold Z}

\define\ca{\Cal A}

\define\cd{\Cal D}
\define\ce{\Cal E}
\define\cf{\Cal F}

\define\ch{\Cal H}

\define\cl{\Cal L}
\define\cm{\Cal M}

\define\cs{\Cal S}
\define\ct{\Cal T}

\define\cv{\Cal V}

\define\fD{\frak D}

\define\fF{\frak F}

\define\tc{\ti c}

\define\tG{\ti G}

\define\tZ{\ti Z}

\define\BBD{BBD}
\define\GR{G}
\define\KA{K}
\define\CRG{L1}
\define\CS{L2}
\define\GF{L3}
\define\CRM{L4}
\define\PCS{L5}
\define\CSD{L6}

\subhead 1\endsubhead
Let $\kk$ be an algebraic closure of a finite field $\FF_q$. Let $G=GL_n(\kk)$. The
group $G(\FF_q)=GL_n(\FF_q)$ can be regarded as the fixed point set of the Frobenius
map $F:G@>>>G,(g_{ij})\m(g_{ij}^q)$. Let $\bbq$ be an algebraic closure of the field
of $l$-adic numbers, where $l$ is a prime number invertible in $\kk$. The characters
of irreducible representations of $G(\FF_q)$ over an algebraically closed field of 
characteristic $0$, which we take to be $\bbq$, have been determined explicitly by 
J.A.Green \cite{\GR}. The theory of character sheaves \cite{\CS} tries to produce 
some geometric objects over $G$ from which the irreducible characters of $G(\FF_q)$
can be deduced for any $q$. This allows us to unify the representation theories of 
$G(\FF_q)$ for various $q$. The geometric objects needed in the theory are provided
by intersection cohomology.

Let $X$ be an algebraic variety over $\kk$, let $X_0$ be a locally closed 
irreducible, smooth subvariety of $X$ and let $\ce$ be a local system over $X_0$ (we
say "local system" instead of "$\bbq$-local system"). Deligne, Goresky and 
MacPherson attach to this datum a canonical object $IC(\bar X_0,\ce)$ (intersection
cohomology complex) in the derived category $\cd(X)$ of $\bbq$-sheaves on $X$; this
is a complex of sheaves which extends $\ce$ to $X$ (by $0$ outside the closure 
$\bar X_0$ of $X_0$) in the most economical possible way so that local Poicar\'e 
duality is satisfied. We say that $IC(\bar X_0,\ce)$ is irreducible if $\ce$ is 
irreducible.

Now take $X=G$ and take $X_0=G_{rs}$ to be the set of regular semisimple elements in
$G$. Let $T$ be the group of diagonal matrices in $G$. For any integer $m\ge 1$ 
invertible in $\kk$ we have an unramified $n!m^n$-fold covering

$\p_m:\{(g,t,xT)\in G_{rs}\T T\T G/T;x\i gx=t^m\}@>>>G_{rs},\qua (g,t,xT)\m g$. 
\nl
An irreducible local system $\ce$ on $G_{rs}$ is said to be admissible if it is a 
direct summand of the local system $\p_{m!}\bbq$ for some $m$ as above. The character 
sheaves on $G$ are the complexes $IC(G,\ce)$ for various admissible local systems 
$\ce$ on $G_{rs}$. 

We show how the irreducible characters of $G(\FF_q)$ can be recovered from
character sheaves on $G$. If $A$ is a character sheaf on $G$ then its inverse image
$F^*A$ under $F$ is again a character sheaf. There are only finitely many $A$ (up to
isomorphism) such that $F^*A$ is isomorphic to $A$. For any such $A$ we choose an 
isomorphism $\ph:F^*A@>\si>>A$ and we form the characteristic function 
$\c_{A,\ph}:G(\FF_q)@>>>\bbq$ whose value at $g$ is the alternating sum of traces of
$\ph$ on the stalks at $g$ of the cohomology sheaves of $A$. Now $\ph$ is unique up
to a non-zero scalar hence $\c_{A,\ph}$ is unique up to a non-zero scalar. It turns
out that

(a) {\it $\c_{A,\ph}$ is (up to a non-zero scalar) the character of an irreducible
representation of $G(\FF_q)$ and $A\m\c_{A,\ph}$ gives a bijection between the set
of (isomorphism classes of) character sheaves on $G$ that are isomorphic to their
inverse image under $F$ and the irreducible characters of $G(\FF_q)$.}
\nl
(This result is essentially contained in \cite{\CRG,\GF}.) The main content of this
result is that the (rather complicated) values of an irreducible character of 
$G(\FF_q)$ are governed by a geometric principle, namely by the procedure which
gives the intersection cohomology extension of a local system.

\subhead 2\endsubhead
More generally, assume that $G$ is a connected reductive algebraic group over $\kk$.
The definition of the $IC(G,\ce)$ given above for $GL_n$ makes sense also in the
general case. The complexes on $G$ obtained in this way form the class of {\it 
uniform} character sheaves on $G$. Consider now a fixed $\FF_q$-rational structure 
on $G$ with Frobenius map $F:G@>>>G$. The analogue of property 1(a) does not hold in
general for $(G,F)$. It is still true that the characteristic functions of the 
uniform character sheaves that are isomorphic to their inverse image under $F$ are 
linearly independent class functions $G(\FF_q)@>>>\bbq$. However they do not form a
basis of the space of class functions. Moreover they are in general not irreducible
characters of $G(\FF_q)$ (up to a scalar); rather, each of them is a linear 
combination with known coefficients of a "small" number of irreducible characters of
$G(\FF_q)$ (where "small" means "bounded independently of $q$"); this result is 
essentially contained in \cite{\CRG,\GF}.

It turns out that the class of uniform character sheaves can be naturally enlarged
to a larger class of complexes on $G$. 

For any parabolic $P$ of $G$, $U_P$ denotes the unipotent radical of $P$. For a 
Borel $B$ in $G$, the images under $c^B:G@>>>G/U_B$ of the double cosets $BwB$ form
a partition $G/U_B=\cup_w(BwB/U_B)$.

An irreducible intersection cohomology complex $A\in\cd(G)$ is said to be a 
character sheaf on $G$ if it is $G$-equivariant and if for some/any Borel $B$ in 
$G$, $c^B_!A$ has the following property: 

($*$) {\it any cohomology sheaf of this complex restricted to any $BwB/U_B$ is a 
local system with finite monodromy of order invertible in $\kk$.}
\nl
Then any uniform character sheaf on $G$ is a character sheaf on $G$. For $G=GL_n$ the 
converse is also true, but for general $G$ this is not so. 

Consider again a fixed $\FF_q$-rational structure on $G$ with Frobenius map 
$F:G@>>>G$. The following partial analogue of property 1(a) holds (under a mild 
restriction on the characteristic of $\kk$).

(a) {\it The characteristic functions of the various character sheaves $A$ on $G$ 
(up to isomorphism) such that $F^*A@>\si>>A$ form a basis of the vector space of 
class functions $G(\FF_q)@>>>\bbq$.}

\subhead 3\endsubhead
We now fix a parabolic $P$ of $G$. For any Borel $B$ of $P$ let
$\tc^B:G/U_P@>>>G/U_B$ be the obvious map. Now $P$ acts on $G/U_P$ by conjugation.

An irreducible intersection cohomology complex $A\in\cd(G/U_P)$ is said to be a 
parabolic character sheaf if it is $P$-equivariant and if for some/any Borel $B$ in
$P$, $\tc^B_!A$ has property 2$(*)$. When $P=G$, we recover the definition of 
character sheaves on $G$.

Consider now a fixed $\FF_q$-rational structure on $G$ with Frobenius map $F:G@>>>G$
such that $P$ is defined over $\FF_q$. Then $G/U_P$ has a natural $\FF_q$-rational
structure with Frobenius map $F$. The following generalization of 2(a) holds (under
a mild restriction on the characteristic of $\kk$).

(a) {\it The characteristic functions of the various parabolic character sheaves $A$
on $G/U_P$ (up to isomorphism) such that $F^*A@>\si>>A$ form a basis of the vector
space $\cv$ of $P(\FF_q)$-invariant functions $G(\FF_q)/U_P(\FF_q)@>>>\bbq$.}
\nl
The proof is given in \cite{\PCS}. It relies on a generalization of property 2(a) to
not necessarily connected reductive groups which will be contained in the series
\cite{\CSD}.

If $h:G(\FF_q)@>>>\bbq$ is the characteristic function of a character sheaf as in
2(a) then by summing $h$ over the fibres of $G(\FF_q)@>>>G(\FF_q)/U_P(\FF_q)$ we
obtain a function $\bar h\in\cv$. It turns out that each function $\bar h$ is a linear
combination of a "small" number of elements in the basis of $\cv$ described above. 
(The fact such a basis of $\cv$ exists is not apriori obvious.)

The parabolic character sheaves on $G/U_P$ are expected to be a necessary ingredient
in establishing the conjectural geometric interpretation of Hecke algebras with
unequal parameters given in \cite{\CRM}.

\subhead 4\endsubhead
In this section $G$ denotes an abelian group with a given family $\fF$ of 
automorphisms such that

(i) if $F\in\fF$ and $n\in\ZZ_{>0}$, then $F^n\in\fF$;

(ii) if $F\in\fF,F'\in\fF$ then there exist $n,n'\in\ZZ_{>0}$ such that
$F^n=F'{}^{n'}$;

(iii) for any $F\in\fF$, the map $G@>>>G,x\m F(x)x\i$ is surjective with finite
kernel.
\nl
For $F\in\fF$ and $n\in\ZZ_{>0}$, the homomorphism 

$N_{F^n/F}:G@>>>G$, $x\m xF(x)\do F^{n-1}(x)$, 
\nl
restricts to a surjective homomorphism $G^{F^n}@>>>G^F$. (If $y\in G^F$ we can find
$z\in G$ with $y=F^n(z)z\i$, by (i),(iii). We set $x=F(z)z\i$. Then $x\in G^{F^n}$ 
and $N_{F^n/F}(x)=y$.) Let $X$ be the set of pairs $(F,\ps)$ where $F\in\fF$ and 
$\ps\in\Hom(G^F,\bbq^*)$. Consider the equivalence relation on $X$ generated by 
$(F,\ps)\si(F^n,\ps\circ N_{F^n/F})$. Let $G^*$ be the set of equivalence classes. 
We define a group structure on $G^*$. We consider two elements of $G^*$; we 
represent them in the form $(F,\ps),(F',\ps')$ where $F=F'$ (using (ii)) and we 
define their product as the equivalence class of $(F,\ps\ps')$; one checks that this
product is independent of the choices. This makes $G^*$ into an abelian group. The 
unit element is the equivalence class of $(F,1)$ for any $F\in\fF$. For $F\in\fF$ we
define an automorphism $F^*:G^*@>>>G^*$ by sending an element of $G^*$ represented 
by $(F^n,\ps)$ with $n\in\ZZ_{>0},\ps\in\Hom(G^{F^n},\bbq^*)$ to $(F^n,\ps\circ F)$
(here $\ps\circ F$ is the composition $G^{F^n}@>F>>G^{F^n}@>\ps>>\bbq^*$); one 
checks that this is well defined. For any $F\in\fF$ the map 
$\Hom(G^F,\bbq^*)@>>>G^*,\ps\m(F,\ps)$ is 

(a) {\it a group isomorphism of $\Hom(G^F,\bbq^*)$ onto the subgroup $(G^*)^{F^*}$ of}
$G^*$.
\nl
(This follows from the surjectivity of $N_{F^n/F}:G^{F^n}@>>>G^F$.)

\subhead 5\endsubhead
Assume now that $G$ is an abelian, connected (affine) algebraic group over $\kk$. We
define the notion of character sheaf on $G$.

Let $\fF$ be the set of Frobenius maps $F:G@>>>G$ for various rational structures on
$G$ over a finite subfield of $\kk$. (These maps are automorphisms of $G$ as an 
abstract group.) Then properties 4(i)-4(iii) are satisfied for $(G,\fF)$ hence the 
abelian group $G^*$ is defined as in \S4. We will give an interpretation of $G^*$ in
terms of local systems on $G$. Let $F\in\fF$. Let $L:G@>>>G$ be the Lang map 
$x\m F(x)x\i$. Consider the local system $E=L_!\bbq$ on $G$. Its stalk at $y\in G$ 
is the vector space $E_y$ consisting of all functions $f:L\i(y)@>>>\bbq$. We have 
$E_y=\op_{\ps\in\Hom(G^F,\bbq^*)}E_y^\ps$ where 

$E_y^\ps=\{f\in E_y;f(zx)=\ps(z)f(x)\qua\frl z\in G^F,x\in L\i(y)\}$. 
\nl
We have a canonical direct sum decomposition $E=\op_\ps E^\ps$ where $E^\ps$ is a
local system of rank $1$ on $G$ whose stalk at $y\in G$ is $E^\ps_y$ ($\ps$ as 
above). There is a unique isomorphism of local systems $\ph:F^*E^\ps@>\si>>E^\ps$ 
which induces identity on the stalk at $1$. This induces for any $y\in G$ the 
isomorphism $E_{F(y)}^\ps@>>>E_y^\ps$ given by $f\m f'$ where $f'(x)=f(F(x))$. If 
$y\in G^F$, this isomorphism is multiplication by $\ps(y)$. Thus, the
characteristic function $\c_{E^\ps,\ph}:G^F@>>>\bbq$ is the character $\ps$. 

Let $n\in\ZZ_{>0}$. Let $L':G@>>>G$ be the map $x\m F^n(x)x\i$. Conider the local
system $E'=L'_!\bbq$ on $G$. Its stalk at $y\in G$ is the vector space $E'_y$ 
consisting of all functions $f':L'{}\i(y)@>>>\bbq$. We define $E_y@>>>E'_y$ by 
$f\m f'$ where $f'(x)=f(N_{F^n,F}x)$ (note that $N_{F^n/F}(L'{}\i(y))\sub L\i(y)$).
This is induced by a morphism of local systems $E@>>>E'$ which restricts to an
isomorphism $E^\ps@>\si>>E'{}^{\ps'}$ where 
$\ps'=\ps\circ N_{F^n/F}\in\Hom(G^{F^n},\bbq^*)$.

From the definitions we see that, if $\ps,\ps'\in \Hom(G^F,\bbq^*)$ then for any 
$y\in G$ we have an isomorphism $E^\ps_y\ot E^{\ps'}_y@>\si>>E^{\ps\ps'}_y$ given by
multiplication of functions on $L\i(y)$. This comes from an isomorphism of local 
systems $E^\ps\ot E^{\ps'}@>\si>>E^{\ps\ps'}$.

A {\it character sheaf} on $G$ is by definition a local system of rank $1$ on $G$ of
the form $E^\ps$ for some $(F,\ps)$ as above. Let $\cs(G)$ be the set of isomorphism
classes of character sheaves on $G$. Then $\cs(G)$ is an abelian group under tensor
product. The arguments above show that $(F,\ps)\m E^\ps$ defines a (surjective) 
group homomorphism $G^*@>>>\cs(G)$. This is in fact an isomorphism. (It is enough to
show that, if $(F,\ps)$ is as above and $\ps'\in\Hom(G^F,\bbq^*)$ is such that the 
local systems $E^\ps,E^{\ps'}$ are isomorphic, then $\ps=\ps'$. As we have seen 
earlier, each of $E^\ps,E^{\ps'}$ has a unique isomorphism $\ph,\ph'$ with its 
inverse image under $F:G@>>>G$ which induces the identity at the stalk at $1$. Then
we must have $\c_{E^\ps,\ph}=\c_{E^{\ps'},\ph'}$ hence $\ps=\ps'$. Note that for 
$F\in\fF$, the map $F^*:G^*@>>>G^*$ corresponds under the isomorphism
$G^*@>\si>>\cs(G)$ to the map $\cs(G)@>>>\cs(G)$ given by inverse image under $F$.
Using this and 4(a), we see that, for $F\in\fF$, the map
$\Hom(G^F,\bbq^*)@>>>\cs(G),\ps\m E^\ps$ is a group isomorphism of
$\Hom(G^F,\bbq^*)$ onto the subgroup of $\cs(G)$ consisting of all character sheaves
on $G$ that are isomorphis to their inverse image under $F$. We see that in this 
case the analogue of 1(a) holds.

From the definitions, we see that, 

(a) {\it if $\cl_1\in\cs(G)$ and $m:G\T G@>>>G$ is the multiplication map then }
$m^*\cl_1=\cl_1\ot\cl_1$.
\nl
In the case where $G=\kk$, our definition of character sheaves on $G$ reduces to that 
of the Artin-Schreier local systems on $\kk$.

\subhead 6\endsubhead
In this section we assume that $G$ is a unipotent algebraic group over $\kk$ of
"exponential type" that is, such that the exponential map from $\Lie G$ to $G$ is 
well defined (and an isomorphism of varieties.) In this case we can define character
sheaves on $G$ using Kirillov theory. Namely, for each $G$-orbit in the dual of 
$\Lie G$ we consider the local system $\bbq$ on that orbit extended by $0$ on the 
complement of the orbit. Taking the Fourier-Deligne transform we obtain (up to 
shift) an irreducible intersection cohomology complex on $\Lie G$ (since the orbit 
is smooth and closed, by Kostant-Rosenlicht). We can view it as an intersection 
cohomology complex on $G$ via the exponential map. The complexes on $G$ thus 
obtained are by definition the character sheaves of $G$. Using Kirillov theory (see
\cite{\KA}) we see that in this case the analogue of 1(a) holds.

Assume, for example, that $G$ is the group of all matrices
$$[a,b,c]=\left(\matrix 1&a&b\\
          0&1&c\\
          0&0&1                               \endmatrix\right)$$
with entries in $\kk$ and that $2\i\in\kk$. Consider the following intersection 
cohomology complexes on $G$:

(i) the complex which on the centre $\{(0,b,0);b\in\kk\}$ is the local system
$\ce\in\cs(\kk),\ce\ne\bbq$ wxtended by $0$ to the whole of $G$;

(ii) the local system $f^*\ce$ where $f[a,b,c]=(a,c)$ and $\ce\in\cs(\kk^2)$.
\nl
The complexes (i),(ii) are the character sheaves of $G$.

\subhead 7\endsubhead
In this section we assume that $G$ is a connected unipotent algebraic group over
$\kk$ (not necessarily of exponential type). We expect that in this case there 
is again a notion of character sheaf on $G$ such that over a finite field, the
characteristic functions of character sheaves form a basis of the space of class 
functions and each characteristic function of a character sheaf is a linear
combination of a "small" number of irreducible characters. Thus here the situation 
should be similar to that for a general connected reductive group rather than that 
for $GL_n$. We illustrate this in one example. Assume that $\kk$ has characteristic
$2$. Let $G$ be the group consisting of all matrices of the form
$$\left(\matrix 1&a&b&c\\
          0&1&d&b+ad\\
          0&0&1&a\\
          0&0&0&1                               \endmatrix\right)$$
with entries in $\kk$; we also write $[a,b,c,d]$ instead of the matrix above. (This
group can be regarded as the unipotent radical of a Borel in $Sp_4(\kk)$.)

Let $\ce_0\in\cs(\kk)$ be the local system on $\kk$ associated in \S5 to $\FF_q$ and
to the homomorphism $\ps_0:\FF_q@>>>\bbq^*$ (composition of the trace  
$\FF_q@>>>\FF_2$ and the unique injective homomorphism $\FF_2@>>>\bbq^*$).

Consider the following intersection cohomology complexes on $G$:

(i) the complex which on the centre $\{[0,b,c,0];(b,c)\in\kk^2\}$ is the local
system $\ce\in\cs(\kk^2),\ce\ne\bbq$ (see \S5) extended by $0$ to the whole of $G$;

(ii) the complex which on $\{[a_0,b,c,0];(b,c)\in\kk^2\}$ (with $a_0\in\kk^*$ fixed)
is the local system $pr_c^*\ce$ where $\ce\in\cs(\kk),\ce\ne\bbq$ (see \S5) extended
by $0$ to the whole of $G$;

(iii) the complex which on $\{[0,b,c,d_0];(b,c)\in\kk^2\}$ (with $d_0\in\kk^*$
fixed) is the local system $f^*\ce_0$ where $f[0,b,c,d_0]=\a b+\a^2d_0c$ (with
$\a\in\kk^*$ fixed) extended by $0$ to the whole of $G$;

(iv) the complex which on $\{[a_0,b,c,d_0];(b,c)\in\kk^2\}$ (with $a_0,d_0\in\kk^*$
fixed) is the local system $f^*\ce_0$ where $f[a_0,b,c,d_0]=a_0^{-2}d_0\i c$
extended by $0$ to the whole of $G$;

(v) the local system $f^*\ce$ on $G$ where $f[a,b,c,d]=(a,d)\in\kk^2$ and
$\ce\in\cs(\kk^2)$.
\nl
By definition, the character sheaves on $G$ are the complexes in (i)-(v) above. Note 
that there are infinitely many subvarieties of $G$ which appear as supports of 
character sheaves (this in contrast with the case of reductive groups). There is a 
symmetry that exchanges the character sheaves of type (ii) with those of type (iii). 
Namely, define $\x:G@>>>G$ by 
$$[a,b,c,d]\m[d,c+ab+a^2d,b^2+dc+abd,a^2].$$
Then $\x$ is a homomorphism whose square is $[a,b,c,d]\m[a^2,b^2,c^2,d^2]$; moreover,
$\x^*$ interchanges the sets (ii) and (iii) and it leaves stable each of the sets 
(i), (iv) and (v). 

Now $G$ has an obvious $\FF_q$-structure with Frobenius map $F:G@>>>G$. We describe
the irreducible characters of $G(\FF_q)$.

(I) We have $q^2$ one dimensional characters $U@>>>\bbq^*$ of the form
$[a,b,c,d]\m\ps_0(xa+yd)$ (one for each $x,y\in\FF_q$).

(II) We have $q-1$ irreducible characters of degree $q$ of the form
$[0,b,c,0]\m q\ps_0(xb)$ (all other elements are mapped to $0$), one for each
$x\in\FF_q-\{0\}$.

(III) We have $q-1$ irreducible characters of degree $q$ of the form
$[0,b,c,0]\m q\ps_0(xc)$ (all other elements are mapped to $0$), one for each
$x\in\FF_q-\{0\}$.

(IV) We have $4(q-1)^2$ irreducible characters of degree $q/2$, one for each
quadruple $(a_0,d_0,\e_1,\e_2)$ where 

$a_0\in\FF_q^*,d_0\in\FF_q^*$, 
$\e_1\in\Hom(\{0,a_0\},\pm 1),\e_2\in\Hom(\{0,d_0\},\pm 1)$, 
\nl
namely

$[a,b,c,d]\m(q/2)\e_1(a)\e_2(d)\ps_0(a_0^{-2}d_0\i(ba+ba_0+c))$,
\nl
if $a\in\{0,a_0\},d\in\{0,d_0\}$; all other elements are sent to $0$.
\nl
A character of type (II) is obtained by inducing from the subgroup
$\{[a,b,c,d]\in G(\FF_q);d=0\}$ the one dimensional character $[a,b,c,0]\m\ps_0(xb)$
where $x\in\FF_q-\{0\}$.
A character of type (III) is obtained by inducing from the commutative subgroup
$\{[a,b,c,d]\in G(\FF_q);a=0\}$ the one dimensional character $[0,b,c,d]\m\ps_0(xc)$
where $x\in\FF_q-\{0\}$.
A character of type (IV) is obtained by inducing from the subgroup
$\{(a,b,c,d)\in G(\FF_q);a\in\{0,a_0\}\}$ (where $a_0\in\FF_q-\{0\}$ is fixed)
the one dimensional character $[a,b,c,d]\m\e_1(a)\ps_0(fd+a_0^{-2}d_0\i(ba+ba_0+c))$
where $f\in\FF_q$ is chosen so that $\ps_0(fd_0)=\e_2(d_0)$ (the induced character
does not depend on the choice of $f$).

Consider the matrix expressing the  characteristic functions of character sheaves
$A$ such that $F^*A\cong A$ (suitably normalized) in terms of irreducible characters
of $G(\FF_q)$. This matrix is square and a direct sum of diagonal
blocks of size $1\T 1$ (with entry $1$) or $4\T 4$ with entries $\pm 1/2$, 
representing the Fourier transform over a two dimensional symplectic $\FF_2$-vector
space. There are $(q-1)^2$ blocks of size $4\T 4$ involving the irreducible 
characters of type IV.

We see that, in our case, the character sheaves have the desired properties. We also
note that in our case, $G(\FF_q)$ has some irreducible character whose degree is not
a power of $q$ (but $q/2$) in contrast with what happens in the situation in \S6.

\subhead 8\endsubhead
Let $\e$ be an indeterminate. For $r\ge 2$ let $\ca_r=\kk[\e]/(\e^r)$. Let 
$G=GL_n(\ca_r)$. Let $B$ (resp.T) be the group of upper triangular (resp. diagonal) 
matrices in $G$. Then $G$ is in a natural way a connected affine algebraic group 
over $\kk$ of dimension $n^2r$ and $B,T$ are closed subgroups of $G$. On $G$ we have
a natural $\FF_q$-structure with Frobenius map $F:G@>>>G$,
$(g_{ij})\m(g_{ij}^{(q)})$ where for $a_0,a_1,\do,a_{r-1}$ in $\kk$ we set
$(a_0+a_1\e+\do+a_{r-1}\e^{r-1})^{(q)}=a_0^q+a_1^q\e+\do+a_{r-1}^q\e^{r-1}$. The 
fixed point set of $F:G@>>>G$ is $GL_n(\FF_q[\e]/(\e^r))$. For $i\ne j$ in $[1,n]$,
we consider the homomorphism $f_{ij}:\kk@>>>T$ which takes $x\in\kk$ to the diagonal
matrix with $ii$-entry equal to $1+\e^{r-1}x$, $jj$-entry equal to $1-\e^{r-1}x$ and
other diagonal entries equal to $1$. Since $T$ is connected and commutative, the 
group $\cs(T)$ is defined (see \S5). Let $\cl\in\cs(T)$. We will assume that $\cl$ 
is {\it regular} in the following sense: for any $i\ne j$ in $[1,n]$, $f_{ij}^*\cl$
is not isomorphic to $\bbq$.

Let $\p:B@>>>T$ be the obvious homomorphism. Consider the diagram
$$G@<a<<Y@>b>>T$$
where 
$$Y=\{(g,xB)\in G\T G/B;x\i gx\in B\}, a(g,xB)=g, b(g,xB)=\p(x\i gx).$$
Then $b^*\cl$ is a local system on $Y$ and we may consider the complex $a_!b^*\cl$ 
on $G$. 

As in \S5, we can find an integer $m_0>0$ such that, for any
$m\in\cm=\{m_0,2m_0,3m_0,\do\}$, $\cl$ is associated to $(\FF_{q^m},\ps_m)$ where 
$\ps_m\in\Hom(T^{F^m},\bbq^*)$. We can regard $\ps_m$ as a character 
$B(\FF_{q^m})@>>>\bbq^*$ via $\p:B@>>>T$; inducing this from $B(\FF_{q^m})$ to 
$G(\FF_{q^m})$ we obtain a representation of $G(\FF_{q^m})$ whose character is 
denoted by $c_m$. It is easy to see (using the regularity of $\cl$) that this 
character is irreducible.

For $m\in\cm$, there is a unique isomorphism $(F^m)^*\cl@>\si>>\cl$ of local systems
on $T$ which induces the identity on the stalk of $\cl$ at $1$. This induces an 
isomorphism $(F^m)^*(b^*\cl)@>\si>>b^*\cl$ (where $F:Y@>>>Y$ is $(g,xB)\m(F(g),F(x)B)$)
and an isomorphism 
$(F^m)^*(a_!b^*\cl)@>\si>>a_!b^*\cl$ in $\cd(G)$. Let $\c_m:G^{F^m}@>>>\bbq$ be the
characteristic function of $a_!b^*\cl$ with respect to this isomorphism. From the 
definitions we see that $\c_m=c_m$. This shows that $a_!b^*\cl$ behaves like a 
character sheaf except for the fact that it is not clear that it is an intersection 
cohomology complex.

We conjecture that:

(a) {\it if $\cl$ is regular then $a_!b^*\cl$ is an intersection cohomology complex
on $G$.}
\nl
(The conjecture also makes sense and is expected to be true when $GL_n$ is replaced
by any reductive group, and $G$ by the corresponding group over $\ca_r$.) Thus one 
can expect that there is a theory of character sheaves for $G$, as far as generic 
principal series representations and their twisted forms is concerned. But one
cannot expect a complete theory of character sheaves in this case (see \S13).

In \S9-\S12 we prove the conjecture in the special case where $G=GL_2(\kk)$ and $r=2$.

\subhead 9\endsubhead
Let $\ca=\ca_2=\kk[\e]/(\e^2)$. Let $V$ be a free $\ca$-module of rank $2$. Let $G$
be the group of automorphisms of the $\ca$-module $V$. This is the group of all 
automorphisms of the $4$-dimensional $\kk$-vector space $V$ that commute with the 
map $\e:V@>>>V$ given by the $\ca$-module structure. Hence $G$ is an algebraic group
of dimension $8$ over $\kk$. Let ${}^0\tG$ be the set of all pairs $(g,V_2)$ where 
$g\in G$ and $V_2$ is a free $\ca$-submodule of $V$ of rank $1$ such that 
$gV_2=V_2$. For $k=1,2$, let $X_k$ be the set of all $\ca$-submodules of $V$ that 
have dimension $k$ as a $\kk$-vector space. Let $\tG$ be the set of all triples 
$(g,V_1,V_2)$ where $g\in G$, $V_1\in X_1,V_2\in X_2$, $V_1\sub V_2$, 
$gV_1=V_1,gV_2=V_2$ and the scalars by which $g$ acts on $V_1$ and $V_2/V_1$ 
coincide. We can regard ${}^0\tG$ as a subset of $\tG$ by $(g,V_2)\m(g,\e V_2,V_2)$.
Note that $\tG$ is naturally an algebraic variety over $\kk$ and ${}^0\tG$ is an 
open subset of $\tG$.

The group of units $\ca'$ of $\ca$ is an algebraic group isomorphic to $\kk^*\T\kk$.
Hence $\cs(\ca')$ is defined. Let $\cl_1\in\cs(\cs'),\cl_2\in\cs(\cs')$. Let 
$\cl=\cl_1\bxt\cl_2\in\cs(\ca'\T\ca')$, $\ce=\cl_2\ot\cl_1^*\in\cs(\ca')$. Define 
$f:{}^0\tG@>>>\ca'\T\ca'$ by $f(g,V_2)=(\a_1,\a_2)$ where $\a_1\in\ca'$ is given by
$gv=\a_1v$ for $v\in V_2$ and $\a_2\in\ca'$ is given by $gv'=\a_2v'$ for 
$v'\in V/V_2$. Let $\tcl=f^*(\cl_1\bxt\cl_2)$, a local system on ${}^0\tG$. Define 
$f_i:{}^0\tG@>>>\ca'$ ($i=1,2$) by $f_1(g,V_2)=\a_1\a_2$, $f_2(g,V_2)=\a_1$ where 
$\a_1,\a_2$ are as above. Then $\tcl=f_1^*\cl_1\ot f_2^*\cl$. (We use 5(a).)

We shall assume that $\cl$ is {\it regular} in the following sense: the restriction
of $\ce$ to the subgroup $\ct=\{1+\e c;c\in\kk\}$ of $\ca'$ is not isomorphic to 
$\bbq$.

\proclaim{Lemma 10} (a) $\tG$ is an irreducible, smooth variety and $\tG-{}^0\tG$ is
a smooth irreducible hypersurface in $\tG$.

(b) We have $IC(\tG,\tcl)|_{\tG-{}^0\tG}=0$.
\endproclaim
Note that $f_1:{}^0\tG@>>>\ca'$ extends to the whole of $\tG$ by 
$f_1(g,V_1,V_2)=\det_{\ca}(g:V@>>>V)$. Hence $f_1^*\cl_1$ extends to a local system
on $\tG$ and we have $IC(\tG,\tcl)=f_1^*\cl_1\ot IC(\tG,f_2^*\ce)$. Hence to prove
(b) it is enough to show that $IC(\tG,f_2^*\ce)$ is zero on $\tG-{}^0\tG$.

Let $Z$ (resp. $H$) be the fibre of the second projection $\tG@>>>X_1$ (resp. 
$\tG-{}^0\tG@>>>X_1$) at $V_1\in X_1$. Since $G$ acts transitively on $X_1$ it is 
enough to show that $Z$ is smooth, irreducible, $H$ is a smooth, irreducible
hypersurface in $Z$ and $IC(Z,f_2^*\ce)$ is zero on $H$ (the restriction of $f_2$ to
$Z$ is denoted again by $f_2$).

Let $e_1,e_2$ be a basis of V such that $V_1=\kk\e e_1$. The subspaces $V_2\in X_2$ 
such that $V_1\sub V_2$ are exactly the subspaces 
$V_2^{z',z''}=\kk\e e_1+\kk(z'e_1+z''\e e_2)$ where $(z',z'')\in\kk^2-\{0\}$. An 
element $g\in G$ is of the form
$$\align&ge_1=a_0e_1+b_0e_2+a_1\e e_1+b_1\e e_2,\\&
ge_2=c_0e_1+d_0e_2+c_1\e e_1+d_1\e e_2\endalign$$
where $a_i,b_i,c_i,d_i\in\kk$ satisfy $a_0d_0-b_0c_0\ne 0$. 

The condition that $g\e e_1\in\kk\e e_1$ is $b_0=0$. 
The condition that $gV_2^{z',z''}=V_2^{z',z''}$ is that $z'b_1+z''d_0=a_0z''$ if 
$z'\ne 0$ (no condition if $z'=0$). 
The condition that the scalars by which $g$ acts on $V_1$ and $V_2^{z',z''}/V_1$ 
coincide is $a_0=d_0$ if $z'=0$ (no condition if $z'\ne 0$). 

We see that we may identify $Z$ with
$$\align\{(a_0,c_0,d_0,a_1,b_1,c_1,d_1;z',z'')&\in\kk^7\T(\kk^2-\{0\})/\kk^*;\\&
a_0\ne 0,d_0\ne 0,z'b_1=z''(a_0-d_0)\}\endalign$$
and $H$ with the subset defined by $z'=0$. In this description it is clear that $Z$
is irreducible, smooth and $H$ is a smooth, irreducible hypersurface in $Z$. The 
function $f_2$ takes a point with $z'\ne 0$ to $a_0+\e(a_1+z''z'{}\i c_0)$. To prove
the statement on intersection cohomology we may replace $Z$ by the open subset 
$z''\ne 0$ containing $H$. Thus we may replace $Z$ by
$$Z_1=
\{(a_0,c_0,d_0,a_1,b_1,c_1,d_1;z)\in\kk^7\T\kk;a_0\ne 0,d_0\ne 0,zb_1=a_0-d_0\}$$
and $H$ by the subset defined by $z=0$. The function $f_2$ is defined on $Z_1-H$ by
$$a_0+\e(a_1+z\i c_0)=(a_0+\e a_1)(1+\e z\i c_0a_0\i).$$
Thus $f_2=f_3f_4$ where $f_3$ (resp. $f_4$) is defined on $Z_1-H$ by $a_0+\e a_1$ 
(resp. $1+\e z\i c_0a_0\i$). Hence $f_2^*\ce=f_3^*\ce\ot f_4^*\ce$. Now $f_3$ 
extends to $Z_1$ hence $f_3^*\ce$ extends to a local system on $Z_1$. We have
$IC(Z_1,f_3^*\ce\ot f_4^*\ce)=f_3^*\ce\ot IC(Z_1,f_4^*\ce)$. It is enough to show 
that $IC(Z_1,f_4^*\ce)$ is zero on $H$. We make the change of variable $c=c_0a_0\i$.
Then $Z_1$ becomes
$$Z_1=\{(a_0,c,a_1,b_1,c_1,d_1;z)\in\kk^7\T\kk;a_0\ne 0,a_0-zb_1\ne 0\},$$
$H$ is the subset defined by $z=0$ and $f_4:Z_1-H@>>>\ca'$ is given by $1+\e z\i c$.
Let $\tZ_1=\{(a_0,c,a_1,b_1,c_1,d_1;z)\in\kk^7\T\kk\}$ and let $H_1$ be the subset 
of $\tZ_1$ defined by $z=0$. Then $Z_1$ is open in $\tZ_1$ and $f_4$ is well defined
on $\tZ_1-H_1$ by $1+\e z\i c$. Hence $f_4^*\ce$ is well defined on $\tZ_1-H_1$. It
is enough to show that $IC(\tZ_1,f_4^*\ce)$ is zero on $H_1$. Let 
$H'=\{(c,z)\in\kk^2;z=0\}$ and define $f':\kk^2-H'@>>>\ca'$ by $f'(c,z)=1+\e z\i c$.
It is enough to show that $IC(\kk^2,f'{}^*\ce)$ is zero on $H'$. Let $P$ be the 
projective line associate to $\kk^2$. Then $H'$ defines a point $x_0\in P$. Since 
$f'$ is constant on lines, it defines a map $h:P-\{x_0\}@>>>\ca'$. Since $P$ is 
$1$-dimensional we have $IC(P,h^*\ce)=\cf$ where $\cf$ is a constructible sheaf on 
$P$ whose restriction to $P-\{x_0\}$ is $h^*\ce$. It is enough to show that 

(c) the stalk of $\cf$ at $x_0$ is $0$;

(d) $H^i(P,\cf)=0$ for $i=0,1$.
\nl
(Indeed, (c) implies that $IC(\kk^2,f'{}^*\ce)$ is zero at $(c,0)$ with $c\ne 0$ and
(d) implies that $IC(\kk^2,f'{}^*\ce)$ is zero at $(0,0)$.)

Consider the standard $\FF_q$-rational structures an $\kk^2,X,\ca'$ and let $F$ be 
the corresponding Frobenius map. We may assume that $\ce$ is associated as in \S5 to
$(\FF_q,\ps)$ where $\ps\in\Hom(\ca'{}^F,\bbq^*)$. For any $m\in\ZZ_{>0}$ there is a
unique isomorphism $\ph_m:(F^m)^*\ce@>\si>>\ce$ which induces the identity on the 
stalk of $\ce$ at $1$. The characteristic function of $\ce$ with respect to this 
isomorphism is $a'\m\ps(N_{F^m/F}(a'))$, $a'\in\ca'{}^{F^m}$. Since, by assumption,
$\ce|_\ct$ is not isomorphic to $\bbq$, $\ps|_{\ct^F}$ is not the trivial character.
Hence $\ps\circ N_{F^m/F}:\ca'{}^{F^m}@>>>\bbq^*$ is non-trivial on $\ct^{F^m}$. Now
$\ph_m$ induces an isomorphism $\ph'_m:(F^m)^*h^*\ce@>\si>>h^*\ce$. We show that

(e) $\sum_{x\in P^{F^m}-\{x_0\}}\tr(\ph'_m,(h^*\ce)_x)=0$.
\nl
An equivalent statement is:

$\sum_{(c,z)\in\FF_{q^m}\T\FF_{q^m}^*}(\ps\circ N_{F^m/F})(1+\e z\i c)=0$,
\nl
which follows from the fact that $\ps\circ N_{F^m/F}:\ca'{}^{F^m}@>>>\bbq^*$ is 
non-trivial on $\ct^{F^m}$. Introducing (e) in the trace formula for Frobenius, we 
see that

(f) $\sum_{i=0}^2(-1)^i\tr(\ph'_m,H^i(P,\cf))=\tr(\ph'_m,\cf_{x_0})$
\nl
where $\cf_{x_0}$ is the talk of $\cf$ at $x_0$ and $\ph'_m$ is in fact equal to 
$\ph'_1{}^m$ (for $m=1,2,3,\do$). By Deligne's purity theorem, $H^i(P,\cf)$ together
with $\ph'_1$ is pure of weight $i$; by Gabber's theorem \cite{\BBD}, $\cf_{x_0}$ 
together with $\ph'_1$ is mixed of weight $\le 0$. Hence from (f) we deduce that 
$H^1(P,\cf)=0,H^2(P,\cf)=0$ and $\dim H^0(P,\cf)=\dim\cf_{x_0}$. By the hard 
Lefschetz theorem \cite{\BBD} we have $\dim H^0(P,\cf)=\dim H^2(P,\cf)$. It follows
that $H^0(P,\cf)=0$ hence $\cf_{x_0}=0$. This proves (c),(d). The lemma is proved.

\proclaim{Lemma 11} Define $\r:{}^0\tG@>>>G$ by $(g,V_2)\m g$. Let $K=\r_!\tcl$. Let
$G_0$ be the open dense subset of $G$ consisting of all $g\in G$ such that
$g:\e V@>>>\e V$ is regular, semisimple. Let $\r_0:\r\i(G_0)@>>>G_0$ be the 
restriction of $\r$. Then $\r_{0!}\tcl$ is a local system on $G_0$. We have
$\dim\supp\ch^iK<\dim G-i$ for any $i>0$.
\endproclaim
The first assertion of the lemma follows from the fact that $\r_0$ is a double 
covering. To prove the second assertion it is enough to show that, for $i>0$, the 
set $G_i$ consisting of the points $g\in G$ such that $\dim\r\i(g)=i$ and 
$\op_jH^j_c(\r\i(g),\tcl)\ne 0$ has codimension $>2i$ in $G$.

Consider the fibre $\r\i(g)$ for $g\in G$. We may assume that, with respect to a 
suitable $\ca$-basis of $V$, $g$ can be represented as an upper triangular matrix 
$\left(\sm a &b\\ 0&c\esm\right)$ with $a,c$ in $\ca'$ and $b\in\ca$. (Otherwise, 
$\r\i(g)$ is empty.) There are five cases:

{\it Case 1.} $a-d\in\ca'$. Then $\r\i(g)$ consists of two points.

{\it Case 2.} $a-d\in\e\ca,b\in\ca'$. Then $\r\i(g)$ is an affine line.

{\it Case 3.} $a-d\in\e\ca-\{0\},b\in\e\ca$. Then $\r\i(g)$ is a disjoint union of
two affine lines.

{\it Case 4.} $a=d,b\in\e\ca-\{0\}$. Then $\r\i(g)$ is an affine line.

{\it Case 5.} $a=d,b=0$. Then $\r\i(g)$ is an affine line bundle over a projective 
line.
\nl
In case 2, we may identify $\r\i(g),\tcl|_{\r\i(g)}$ with 
$P-\{x_0\},\cf|_{P-\{x_0\}}$ in the proof of Lemma 10. Then the argument in that 
proof shows that $H^j_c(\r\i(g),\tcl)=0$ for all $j$. We see that $G_1$ consists of
all $g$ as in case 3 and 4, hence $G_1$ has codimension $3$ in $G$. We see that 
$G_2$ consists of all $g$ as in case 5, hence $G_2$ has codimension $6$ in $G$. The
lemma is proved. Note that without the assumption that $\cl$ is regular, the last 
assertion of the lemma would not hold (there would be a violation coming from $g$ in
case 2.)

\subhead 12\endsubhead
We show:
$$\r_!\tcl=IC(G,\r_{0!}\tcl).\tag a$$
Define $\ti\r:\tG@>>>G$ by $\ti\r(g,V_1,V_2)=g$. Clearly, $\ti\r$ is proper. Let 
$j:{}^0\tG@>>>G$ be the inclusion. We have $\r=\ti\r\circ j$ hence
$\r_!\tcl=\ti\r_!(j_!\tcl)$. By Lemma 10, we have $j_!\tcl=IC(\tG,\tcl)$ hence
$\r_!\tcl=\ti\r_!IC(\tG,\tcl)$. Since $\ti\r$ is proper, $\ti\r_!$ commutes with the
Verdier duality $\fD$. Hence $\fD(\r_!\tcl)=\ti\r_!\fD IC(\tG,\tcl)$. Hence
$\fD(\r_!\tcl)$ equals $\ti\r_!IC(\tG,\tcl^*)$ up to a shift. Now the same 
argument that shows $j_!\tcl=IC(\tG,\tcl)$ shows also $j_!\tcl^*=IC(\tG,\tcl^*)$. 
Hence, up to shift, $\fD(\r_!\tcl)$ equals $\ti\r_!j_!\tcl^*=\r_!\tcl^*$. Now 
the argument in Lemma 12 can also be applied to $\tcl^*$ instead of $\tcl$ and 
yields $\dim\supp\ch^i\r_!\tcl^*<\dim G-i$ for any $i>0$. Thus, $\r_!\tcl$
satisfies the defining properties of $IC(G,\r_{0!}\tcl)$ hence it is equal to it. 
This proves (a).

We see that conjecture 8(a) holds for $n=2,r=2$.

\subhead 13\endsubhead
If $G$ is a connected affine algebraic group over $\kk$ which is neither
reductive nor nilpotent, one cannot expect to have a complete theory character sheaves 
for $G$. Assume for example that $G$ is the group of all matrices
$$[a,b]=\left(\sm a&b\\
          0&1                               \esm\right)$$
with entries in $\kk$. The group $G(\FF_q)$ (for the obvious $\FF_q$-rational
structure) has $(q-1)$ one dimensional representations and one $(q-1)$-dimensional
irreducible representation. The character of a one dimensional representation can be
realized in terms of an intersection cohomology complex (a local system on $G$), but 
that of the $(q-1)$ dimensional irreducible representation appears as a difference of 
two intersection cohomology complexes, one given by the local system $\bbq$ on the 
unipotent radical of $G$ and one supported by the unit element of $G$. A similar 
phenomenon occurs for $G$ as in \S9 and for a $(q^2-1)$-dimensional irreducible 
representation of $G(\FF_q)$.

\enddocument